\documentclass[twoside]{amsart}
\usepackage{latexsym}
\usepackage{amssymb,amsmath,amsopn}
\usepackage[dvips]{graphicx}   
\usepackage{color,epsfig}      

               {\begin{list}{}{\leftmargin#1\rightmargin#2}\item{}}%
               {\end{list}}

\theoremstyle{definition}

\def\bea{\begin{eqnarray}}
\def\eea{\end{eqnarray}}

\usepackage[normalem]{ulem}

\begin{document}
\title[Building a monostable tetrahedron]{Building a monostable tetrahedron}

\author[G. Alm\'adi, R. J. MacG. Dawson and G. Domokos] {Gerg\H o Alm\'adi, Robert J. MacG. Dawson and G\'abor Domokos}
\address{Gerg\H o Alm\'adi, HUN-REN-BME Morphodynamics Research Group, Budapest University of Technology and Economics,
M\H uegyetem rakpart 1-3., Budapest, Hungary, 1111}
\email{gergo.almadi14@gmail.com}
\address{Robert J. MacG. Dawson, Dept. Mathematics and Computer Science, Saint Mary's University, Halifax, Nova Scotia B3H 3C3, Canada}
\email{rdawson@cs.stmarys.ca}
\address{G\'abor Domokos, Dept. of Morphology and Geometric Modeling and HUN-REN-BME Morphodynamics Research Group, Budapest University of Technology and Economics,
M\H uegyetem rakpart 1-3., Budapest, Hungary, 1111}
\email{domokos@iit.bme.hu}

\thanks{GD and GA: Support of the NKFIH Hungarian Research Fund grant 149429 and of grant BME FIKP-V\'IZ by EMMI is kindly acknowledged. AG:  This research has been supported by the
program UNKP-23-3 by ITM and NKFIH. The gift representing the Albrecht Science Fellowship
is gratefully appreciated.}

\begin{abstract}
In this short note we describe what we believe to be the first working model of a monostable tetrahedron. 
\end{abstract}
\maketitle

\section{Background and motivation}\label{sec:mot}
In the SIAM Review problems column for July 1966 J.H. Conway and R. Guy asked
readers \cite{Conway1} to show that a homogeneous tetrahedron, of any shape whatsoever,
must be gravitationally stable on at least two faces; this was answered
in the same column in 1969 \cite{Conway2}. In 1967 Heppes provided an example
of a homogeneous tetrahedron which would tip over two of its edges before coming to
rest on a stable face \cite{Heppes1967}. In 1984, Conway told the second author,
that he had shown, some time before, that it was possible for a suitably-shaped \emph{nonhomogeneous}
tetrahedron to be stable on only one face. While the proof did not appear
to be particularly difficult, at that time it was not formalized. 

What did appear as a challenge, though, was a physical realization of such an object. 
The second author built a model (now lost) from lead foil and finely-split bamboo, which appeared
to tumble sequentially from one face, through two others, to its final
resting position. However, the model was suspect: the edges bent under
the weight of the lead pellet. In this problem \emph{bending the edges} is equivalent
to \emph{bending the rules}: if we are permitted to use curved edges then
it is very easy to build a monostable tetrahedron.  Consider a tetrahedron
the edges of which are running on a sphere. If the center of mass is not identical
to the center of the sphere then the tetrahedron will be monostable.
(In fact, this holds for any polyhedron.)

\section{The geometry of monostable tetrahedra}

The theoretical aspects of the problem were substantially advanced 
in  \cite{Almadi}, where not only the formal proof for the existence
of a monostable tetrahedron was delivered but many other interesting
properties of such objects were described. 

If a tetrahedron can have a mass distribution which makes it monostable then we call it \emph{loadable}. In \cite{Almadi} it was shown that loadable tetrahedra are characterized by the existence of an \emph{obtuse path}, consisting of three obtuse edges $\overline{ab}$,$\overline{bc}$,$\overline{cd}$. Let $A$ be the face opposite $a$, etc.: if an obtuse path exists, then it is always possible to achieve any of the \emph{falling patterns} $C\rightarrow D \rightarrow A \rightarrow B$, $C\rightarrow D \rightarrow A \leftarrow B$, $C\rightarrow D \leftarrow A \leftarrow B$, or $C\leftarrow D \leftarrow A \leftarrow B$ by appropriate choice of the center of mass.  The falling patterns ending on $A$ and $D$ are called ``Type I'', the others ``Type II.'' The set of centers of mass leading to a given falling pattern is called its \emph{loading zone}; this is always a tetrahedron.

\section{Building the model}
Using all these geometrical facts we built, what we believe is the first
working physical model of a monostable tetrahedron.

\begin{figure}[h!]
\centering
    \includegraphics[width=0.95\textwidth]{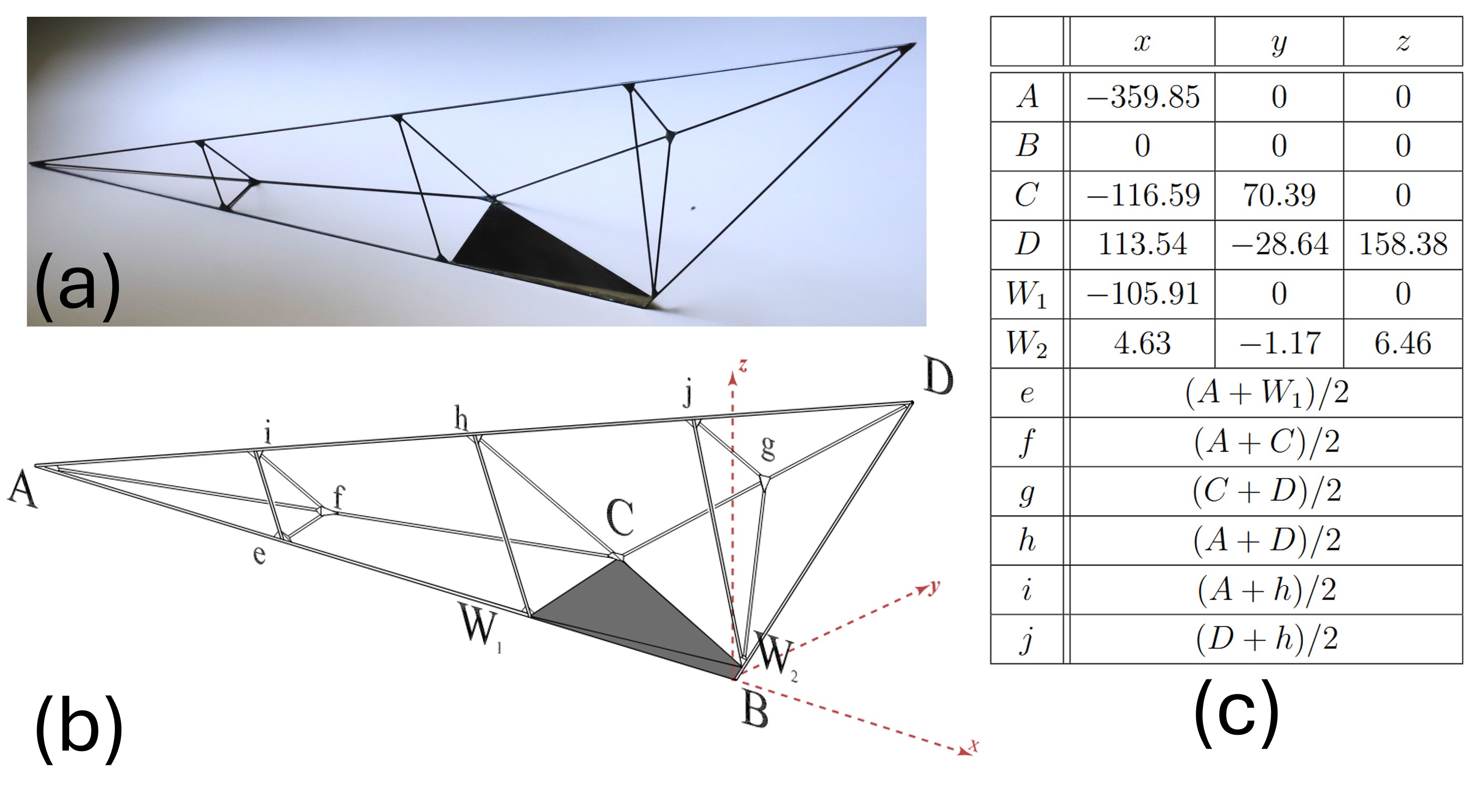}
\caption{The geometry of a working monostable tetrahedron construction. (a) Photo of the structure. Frame made of pultruded carbon fibre tubes  with outer diameter = 1mm, inner diameter = 0,5mm, density 1,36 g/cm$^3$. Core made of tungsten-carbide of density 14,15 g/cm$^3$. Edges are joined with an epoxy glue of density 1,3g/cm$^3$. (b) Axonometric drawing. (c) Orthogonal $[xyz]$ coordinates of nodes in expressed in mm.  }
\label{fig:1}
\end{figure}

Numerical evidence (see Table \ref{tab:1}) indicated that the  volume of Type I loading zones
exceeded the size of Type II loading zones by at least one order
of magnitude, so we aimed for a Type I model.

\begin{table}\begin{center}
\begin{tabular}{|l|c||c|}
\hline
Falling pattern & Type & Volume of loading zone  \\
\hline
\hline
$B \rightarrow A \rightarrow D \leftarrow C$ &  I. & 1.4318 cm$^3$ \\
\hline
$C \rightarrow D \rightarrow A \leftarrow B$  & I. & 0.5716 cm$^3$ \\
\hline
$B \rightarrow A \rightarrow D \rightarrow C$  & II. & 0.0199 cm$^3$  \\
\hline
$C \rightarrow D \rightarrow A \rightarrow B$  &  II. & 0.0067 cm$^3$\\
\hline
\end{tabular}\caption{Falling patterns and loading zone volumes for the tetrahedron in Figure \ref{fig:1}.
Volume of entire tetrahedron is 668.624 cm$^3$}\label{tab:1}
\end{center}\end{table}

To place the center of mass into the loading zone, we adopted
a straightforward geometrical approach: we imagined that the tetrahedron
is made from two different materials: one heavy and the other light.
We assumed that the interface between the two materials is planar
and we used the geometry of the loading zones to design the
planar interface. 
Initial computations showed that with the chosen planar interface, if the density
of the heavy part exceeded the density of the light part
by at least three orders of magnitude then the model still
appeared functional.

To achieve this high density ratio we constructed the tetrahedron
as a skeleton of very light carbon-fibre tubes, so the
density of the light part could be derived from the mass
of the skeleton divided by the volume of the tetrahedron.
The heavy part was constructed from tungsten carbide ($14.15 \mbox{g}/\mbox{cm}^3$, denser than
lead). Using these density values we could compute the minimal size
at which the model still appeared to be functional and in the design we selected 10\%
larger linear size.

The resulting object (illustrated in Figure \ref{fig:1}) is stable only on one of its faces,
on face $D$ (opposite vertex $D$). If
placed on any other face, it will tumble to $D$ according to the
directed graph $A \rightarrow B \rightarrow D \leftarrow C$.
These falling patterns are shown in \cite{video}.

\section{Open questions and potential applications}
We also made calculation for the Type II falling pattern.
In table \ref{tab:1} we summarized the volumes
of the 4 loading zones associated with the tetrahedron in Figure \ref{fig:1}.
We found that considering the carbon tube frame with identical geometry,
for the pattern $B \rightarrow A \rightarrow D \rightarrow C$  the
weight would have to have a density of at least $234
\mbox{g}/\mbox{cm}^3$, an order of magnitude denser than any material
found on Earth. This shows that the hand-made split bamboo model,
displaying the Type II pattern, must have relied on the curved geometry of its edges.

\begin{figure}
    \centering
    \includegraphics[width=0.7\linewidth]{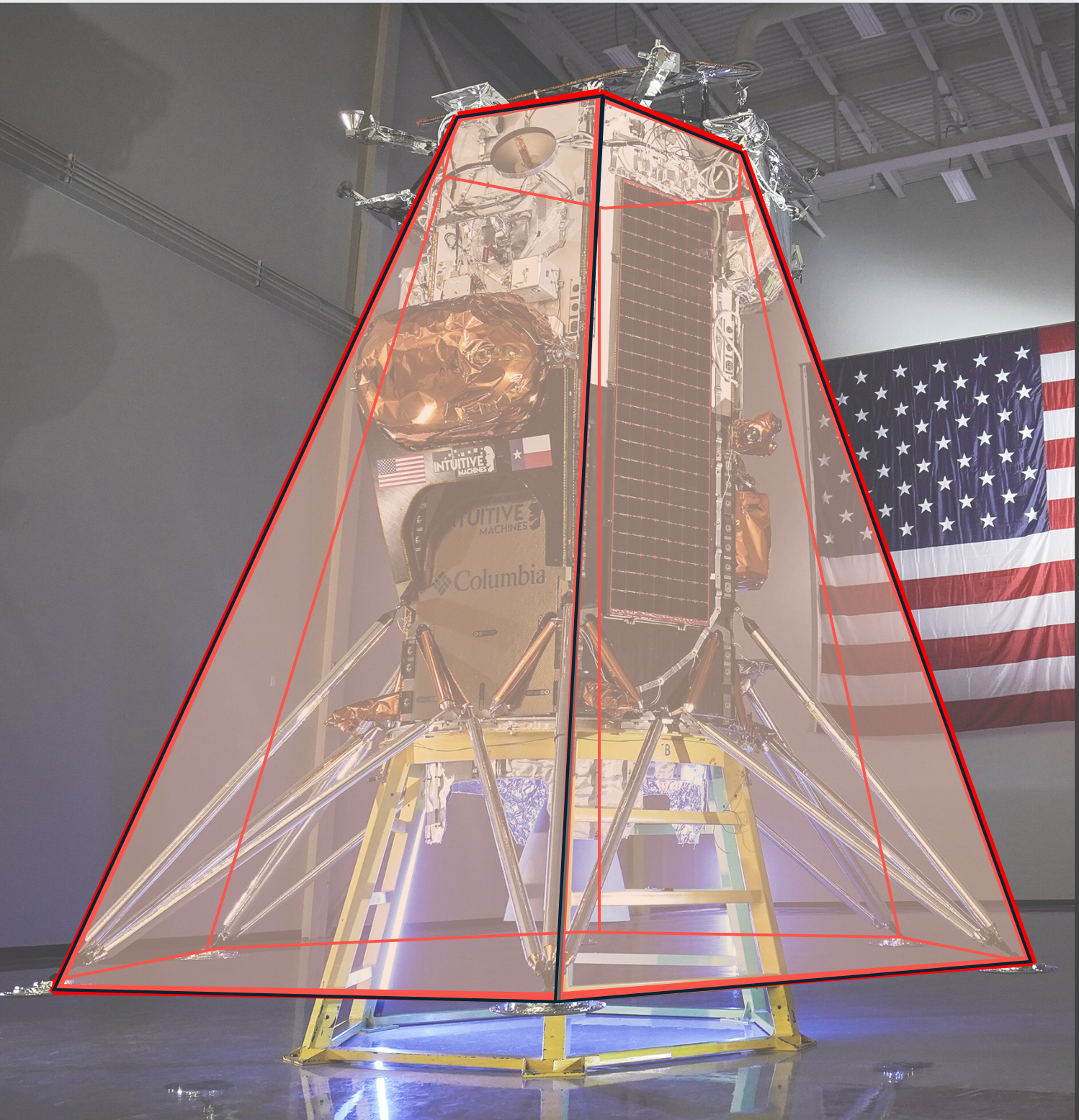}
    \caption{The approximation of the convex hull of the lunar lander Nova-C \emph{Odysseus} as an inhomogeneous convex polyhedron. In case of a failure, like on mission \emph{IM-1} where \emph{Odysseus} tipped over \cite{IM1}, the lander rolls on its convex hull.}
    \label{fig:2}
\end{figure}

As pointed out in Section \ref{sec:mot}, in this problem
it is important to consider straight edges.
At first sight, the latter may appear as a mathematical abstraction, nevertheless,
they play a central role in many applications. If we consider
the stability of a non-convex object supported on a horizontal plane then
equilibrium positions will be associated with its convex hull which,
for many non-smooth objects, is a polyhedron with ideal (straight) edges.
So, studying the balancing properties of inhomogeneous polyhedra with ideal
(straight) edges is a real problem in engineering.

The latter type of problem was recently shown to be relevant
when three lunar landers tipped over \cite{IM1, IM2, SLIM}.
These machines are non-convex. non-smooth objects which would
roll on their respective convex hulls, which, in turn, can be approximated
by convex polyhedra. Figure \ref{fig:2} illustrates
the Odysseus lander from mission IM-1 \cite{IM1} and its convex hull.

While it may not be possible to design objects which
can passively self-right on any terrain, designing
for self-righting on a horizontal support may be feasible
and we hope that for those designs our study could offer insights.

\section*{Acknowledgement}
GA and GD: This research was supported by NKFIH grants K149429 and EMMI FIKP grant VIZ. GA: This research has been supported by the program EK\"OP-24-2, funded by ITM and NKFI.  The gift representing the Albrecht Science Fellowship is gratefully appreciated.

\end{document}